\documentclass[12pt]{article}
 
\usepackage{amssymb,amsmath}
 
\usepackage{accents} 

\usepackage{dsfont}

\usepackage{tikz}
\usetikzlibrary{automata,topaths}
\usetikzlibrary{shapes}
\usetikzlibrary{plotmarks}
 
\usepackage{wasysym} 

\usepackage{color} 
\usepackage[colorlinks=true, urlcolor=blue,linkcolor=blue, citecolor=blue]{hyperref}

\DeclareMathAlphabet{\eufrak}{U}{}{}{} 
\SetMathAlphabet\eufrak{normal}{U}{euf}{m}{n}
\SetMathAlphabet\eufrak{bold}{U}{euf}{b}{n}

\numberwithin{equation}{section}

\newenvironment{Proof}{\removelastskip\par\medskip
\noindent{\em Proof.} \rm}{\penalty-20\null\hfill$\square$\par\medbreak}

\allowdisplaybreaks

 \def\real{{\mathord{\mathbb R}}}
 
 \def\inte{{\mathord{\mathbb N}}}
 
 \def\qu{{\mathord{\mathbb Z}}}

 \def\Var{{\mathrm{{\rm Var}}}}

 \def\real{{\mathord{{\rm I\kern-3pt R}}}}        
 
 \def\inte{{\mathord{{\rm I\kern-3pt N}}}}
 \def\sZZ{{\rm Z\kern-.45em{}Z}}

 \def\sQQ{{\kern 0.27em \vrule height1.45ex width0.03em depth0em
           \kern-0.30em \rm Q}}
 \def\qu{{\mathchoice
         {\sQQ}
         {\sQQ}
   {\kern 0.225em \vrule height1.05ex width0.025em depth0em \kern-0.25em \rm Q}
   {\kern 0.180em \vrule height0.78ex width0.020em depth0em \kern-0.20em \rm Q}
         }}
 \def\sGG{{\kern 0.27em \vrule height1.45ex width0.03em depth0em
           \kern-0.30em \rm G}}
 \def\gg{{\mathchoice
         {\sGG}
         {\sGG}
   {\kern 0.225em \vrule height1.05ex width0.025em depth0em \kern-0.25em \rm G}
   {\kern 0.180em \vrule height0.78ex width0.020em depth0em \kern-0.20em \rm G}
         }}

 \newtheorem{prop}{Proposition}[section]
 \newtheorem{lemma}[prop]{Lemma}
 
 \newtheorem{corollary}[prop]{Corollary}

\def\E{\mathop{\hbox{\rm I\kern-0.20em E}}\nolimits}

 \newcounter{hyp}
\newcommand{\re}{\mathrm{e}}            

\title{\huge Moments of $k$-hop counts in the random-connection model} 
 
\author{
\large 
 Nicolas Privault\thanks{nprivault@ntu.edu.sg} 
\small
Division of Mathematical Sciences 
\\ 
\small 
School of Physical and Mathematical Sciences 
\\ 
\small
Nanyang Technological University 
\\ 
\small 
21 Nanyang Link 
\\ 
\small
Singapore 637371
}

\begin{document}

\hyphenation{func-tio-nals} 
\hyphenation{Privault} 

\maketitle 

\vspace{-0.9cm}

\baselineskip0.6cm
 
\begin{abstract} 
  We derive moment identities for the stochastic integrals of
  multiparameter processes in a random-connection model 
  based on a point process admitting a Papangelou
  intensity. 
  Those identities are written using sums over partitions,
  and they reduce to sums over non-flat partition diagrams
  in case the multiparameter processes vanish on diagonals. 
  As an application, we obtain general identities 
  for the moments of $k$-hop counts in the random-connection model,
  which simplify the derivations available in the literature. 
\end{abstract} 
 
 
\baselineskip0.7cm
 
\section{Introduction} 
The random-connection model, see e.g. Chapter~6 of
\cite{meester}, is a classical model in
continuum percolation.
It consists in a random graph built on 
the vertices of a point process on $\real^d$,
by adding edges between two distinct
vertices $x$ and $y$ with probability $H (\Vert x - y \Vert )$. 
In the case of the Rayleigh fading
$H_\beta (\Vert x - y \Vert ) = \re^{- \beta \Vert x - y \Vert }$
with $x, y \in \real^2$, 
the mean value of the number $N^{x,y}_k$ of $k$-hop paths
connecting $x\in \real^d$ to $y\in \real^d$
has been computed in  \cite{kartun-giles}, together with
the variance of $3$-hop counts. 
However, this argument does not extend to $k\geq 3$
as the proof of the variance identity for $3$-hop counts 
in \cite{kartun-giles} relies on the known Poisson
distribution of the $2$-hop count. 
As shown by \cite{kartun-giles}, 
the knowledge of moments can provide accurate numerical estimates of the
probability $P( N^{x,y}_k >0)$ of at least one $k$-hop path, 
by expressing it as a series of factorial moments, and 
the need for a general theory of such moments has been pointed out
therein. 

\bigskip
 
On the other hand, moment identities for Poisson stochastic integrals with
random integrands have been obtained in 
\cite{momentpoi} based on moment identities for the Skorohod 
integral on the Poisson space, see \cite{priinvcr,prinv},
and also \cite{prob} for a review. 
These moment identities have been extended to
point processes with Papangelou intensities by \cite{flint},
and to multiparameter processes by \cite{bogdan}. 
Factorial moments have also been computed by \cite{bretonprivaultfact}
for point processes with Papangelou intensities.

\bigskip

In this paper we derive closed-form expressions for the moments of
the number of $k$-hop paths in the random-connection model.
In Proposition~\ref{djklssa} 
the moment of order $n$ of the $k$-hop count is given as 
a sum over non-flat partitions of $\{1,\ldots , nk\}$
in a general random-connection model based on a point process admitting
a Papangelou intensity.
Those results are then specialized to the case of Poisson point processes,
with an expression for the variance of the $k$-hop count
given in Corollary~\ref{dvcsdf} using a sum over integer sequences.
Finally, in the case of Rayleigh fadings we show that some results
of \cite{kartun-giles}, 
such as the computation of variance for $3$-hop counts,
can be recovered via a shorter argument, see Corollary~\ref{fdfdskl}. 
  
\bigskip
 
We proceed as follows.
After presenting some background notation on
point processes and Campbell measures, see \cite{kallenberg},
in Section~\ref{s2} we review the derivation
of moment identities for stochastic integrals using sums over partitions.
In the multiparameter case we rewrite those identities for processes
vanishing on diagonals, based on non-flat partition diagrams.
In Section~\ref{s4} we apply those results to the computation of the 
moments of $k$-hop counts in the random-connection model,
and we specialize such computations to the case of Poisson
point processes in Section~\ref{s4.1}. 
Section~\ref{s5} is devoted to explicit computations in the case
of Rayleigh fadings, which result into simpler derivations in
comparison with the current literature on moments
in the random-connection model. 

\subsubsection*{Notation on point processes}

Let $X$ be a Polish space with Borel $\sigma$-algebra ${\cal B}(X)$,
equipped with a $\sigma$-finite non-atomic measure $\lambda (dx)$. 
 We let
 $$
 \Omega^X = \big\{
 \omega = \{ x_i \}_{i\in I} \subset X \ : \
 \#( A \cap \omega ) < \infty 
 \mbox{ for all compact } A\in {\cal B} (X) 
 \big\}
 $$
 denote the space of locally finite configurations on $X$, whose elements 
 $\omega \in \Omega^X$ are identified with the Radon point measures 
 $\displaystyle \omega = \sum_{x\in \omega} \epsilon_x$, 
 where $\epsilon_x$ denotes the Dirac measure at $x\in X$. 
 A point process is a probability measure $P$ on $\Omega^X$ equipped 
 with the $\sigma$-algebra ${\cal F}$ generated by the topology of 
 vague convergence. 

 \bigskip 
 
 Point processes can be 
 characterized by their Campbell measure $C$ defined on 
 ${\cal B}(X) \otimes {\cal F}$ by 
$$
C(A\times B) : =\Bbb{E}\left[\int_A \mathds{1}_B(\omega\setminus\{x\})\
\omega(dx)\right], 
 \quad 
 A\in{\cal B}(X), \quad B\in{\cal F}, 
$$ 
 which satisfies the Georgii-Nguyen-Zessin \cite{nguyen} identity 
 
\begin{equation}
\label{eq:GNZ0}
 \Bbb{E} \left[ 
 \int_X u ( x ; \omega ) 
 \omega ( dx ) 
 \right] 
 = 
 \int_{\Omega^X} 
 \int_X u ( x; \omega \cup x ) 
C(dx,d\omega ) 
, 
\end{equation}
for all measurable processes
$u: X\times\Omega^X\to \real$ such that both sides of \eqref{eq:GNZ0} make sense. 

 \bigskip 
 
 In the sequel we deal with point processes whose Campbell 
 measure $C(dx,d \omega )$ is absolutely continuous with respect to 
 $\lambda \otimes P$, i.e. 
$$ 
 C(dx,d\omega ) = c(x;\omega ) \lambda (dx ) P(d\omega ), 
$$ 
 where the density $c(x;\omega)$ is called the Papangelou density. 
 We will also use the random measure 
 $\hat{\lambda}^n ( d\eufrak{x}_n )$ defined on $X^n$ by 
$$ 
 \hat{\lambda}^n ( d\eufrak{x}_n ) 
 = 
 \hat{c} ( \eufrak{x}_n ;\omega ) 
 \lambda (dx_1) \cdots \lambda (dx_n) 
, 
$$ 
 where $\hat{c} ( \eufrak{x}_n ;\omega )$
 is the compound Campbell density 
$
\hat{c}: \Omega_0^X\times\Omega^X \longrightarrow \real_+ 
$
defined inductively
on the set $\Omega_0^X$ of finite configurations 
in $\Omega^X$ by
\begin{equation}
\label{eq:hatc}
 \hat{c} ( 
 \{x_1,\ldots ,x_n , y \} ; \omega ) 
 : = 
 c( y ; \omega ) 
 \hat{c} ( 
 \{x_1,\ldots ,x_n \} ; \omega \cup \{ y\} ) 
, \qquad n \geq 0,  
\end{equation} 
see Relation~(1) in \cite{flint}.
In particular, the Poisson point process with intensity $\lambda (dx)$ 
 is a point process with Campbell measure $C=\lambda \otimes P$ 
 and $c(x; \omega)=1$, and in this case the identity \eqref{eq:GNZ0}
 becomes the Slivnyak-Mecke formula \cite{slivnyak}, \cite{jmecke}. 
Determinantal point processes are examples of 
point processes with Papangelou intensities, see e.g.
Theorem~2.6 in \cite{dfpt2}, and
they can be used for the modeling
of wireless networks with repulsion,
see e.g. \cite{haenggi}, \cite{MSa}, \cite{flpnw4}. 

\section{Moment identities} 
\label{s2}
 The moment of order $n\geq 1$ of
 a Poisson random variable $Z_\alpha$ with parameter $\alpha > 0$
 is given by 
\begin{equation} 
\label{ws3} 
 \E [ Z^n_\alpha ] 
 =
 \sum_{k=0}^n \alpha^k S(n,k), 
\qquad n \in \inte, 
\end{equation} 
 where $S(n,k)$ is the number of ways to partition a set of $n$ objects 
 into $k$ non-empty subsets, 
 see e.g. Proposition~3.1 of \cite{boyadzhiev}. 
 Regarding Poisson stochastic integrals of
 deterministic integrands, in \cite{bassan} the moment formula 
\begin{equation} 
\label{**} 
 \E \left[ 
 \left( \int_X f (x) \omega (dx) \right)^n 
 \right] 
 = 
 n! 
 \! \! \! \! \! 
 \sum_{ 
  r_1+2r_2+\cdots + n r_n =n 
 \atop 
 r_1,\ldots , r_n \geq 0
} 
 \prod_{k=1}^n 
 \left( 
 \frac{1}{(k!)^{r_k}r_k!} 
 \left( \int_X f^k (x) \lambda (dx) \right)^{r_k} 
 \right) 
\end{equation} 
has been proved for deterministic
functions $f \in \bigcap_{p\geq 1} L^p(X,\lambda )$. 

 \bigskip

 The identity \eqref{**} has been rewritten in the langage of sums over partitions, 
 and extended to Poisson stochastic integrals of random integrands in Proposition~3.1
 of \cite{momentpoi}, 
 and further extended to point processes admitting a Panpangelou intensity
 in  Theorem~3.1 of \cite{flint}, see also \cite{bretonprivaultfact}. 
 In the sequel, given $\eufrak{z}_n =(z_1, \ldots, z_n)\in X^n$, we will use 
 the shorthand notation $\varepsilon^+_{\eufrak{z}_n }$ 
 for the operator 
$$ 
 ( \varepsilon^+_{\eufrak{z}_n } F)(\omega)=F(\omega\cup\{z_1, \ldots, z_n\}), 
 \qquad 
 \omega \in \Omega, 
$$ 
 where $F$ is any random variable on $\Omega^X$.
 Given $\rho = \{ \rho_1,\ldots, \rho_k \} \in \Pi [n]$
 a partition of $\{1,\ldots , n\}$ of size $|\rho|=k$, we let $|\rho_i |$
 denote the cardinality of each block $\rho_i$, $i=1,\ldots , k$. 
\begin{prop} 
\label{djld} 
 Let $u : X \times \Omega^X \longrightarrow \real$ be 
 a (measurable) process. For all $n\geq 1$ we have 
\begin{equation} 
\nonumber 
 \Bbb{E} \left[ 
 \left( \int_X u(x ; \omega ) \omega  (dx) \right)^n 
 \right] 
 =
 \sum_{ \rho \in \Pi [n] }
 \Bbb{E} \left[ 
   \int_{X^{|\rho |}} \epsilon^+_{\eufrak{z}_{|\rho|} } \prod_{l=1}^{|\rho |} u^{|\rho_l|} (z_l )
   \hat{\lambda}^{|\rho |} ( d \eufrak{z}_{|\rho|} ) \right], 
\end{equation} 
 where the sum runs over all partitions 
 $\rho$ of $\{ 1 , \ldots , n \}$ with cardinality $| \rho |$.
\end{prop} 
Proposition~\ref{djld} has also been extended, 
together with joint moment identities, 
to multiparameter processes
$(u_{z_1,\ldots z_r})_{(z_1,\ldots z_r)\in X^r}$, 
see Theorem~3.1 of \cite{bogdan}.
For this, let $\Pi [n\times r]$ denote the set of all partitions
of the set
 $$
 \Delta_{n\times r} := \{1,\ldots , n\} \times \{1,\ldots , r\}
 =
 \big\{ (k,l) \ : \
 k=1,\ldots , n, \ l = 1,\ldots , r \big\}, 
 $$
  identified to $\{1,\ldots , nr\}$, and let
$\pi : = (\pi_1,\ldots , \pi_n) \in \Pi [n\times r]$ denote the
 partition made of the $n$ blocks
$\pi_k := \{ (k,1), \ldots , (k,r) \}$ of size $r$, for 
$k=1,\ldots , n$.
Given $\rho = \{ \rho_1,\ldots , \rho_m \}$ a partition of
$\Delta_{n\times r}$, we let $\zeta^\rho : \Delta_{n\times r} \longrightarrow
\{1,\ldots , m \}$ denote the mapping defined as
$$
\zeta^\rho (k,l)=p \mbox{ ~if and only if } 
(k,l)\in \rho_p,
$$
$k=1,\ldots , n$, $l=1,\ldots , r$,
$p=1,\ldots , m$.
In other words, $\zeta^\rho (k,l)$
denotes the index $p$ of the block
$\rho_p\subset \Delta_{n\times r}$ to which $(k,l)$ belongs.
\\
 
Next, we restate Theorem~3.1 of \cite{bogdan}
by noting that, in the same way as in Proposition~\ref{djld}, 
it extends to point processes admitting a Papangelou intensity
using the arguments of \cite{flint}, \cite{bretonprivaultfact}. 
 When 
 $(u(z_1,\ldots ,z_k ; \omega ) )_{z_1,\ldots , z_k \in X}$  
 is a multiparameter process, we will write 
$$
\epsilon^+_{\eufrak{z}_k} u(z_1,\ldots ,z_k ; \omega ) 
:=
u\big(
z_1,\ldots ,z_k ; \omega \cup \{
 z_1,\ldots , z_k \}
 \big),
 \quad
 \eufrak{z}_n =(z_1, \ldots, z_n)\in X^n, 
$$
 and in this case we may drop the variable $\omega \in \Omega^X$  
 by writing $\epsilon^+_{\eufrak{z}_k} u(z_1,\ldots ,z_k ; \omega )$
 instead of $\epsilon^+_{\eufrak{z}_k} u(z_1,\ldots ,z_k ; \omega )$.
\begin{prop} 
  \label{p2.2}
   Let $u : X^r \times \Omega^X \longrightarrow \real$ be 
 a (measurable) $r$-process. We have 
\begin{equation} 
\nonumber 
  \Bbb{E} \left[ 
   \left( \int_{X^r} u(z_1,\ldots , z_r ; \omega ) \omega  (dz_1) \cdots
   \omega (dz_r) \right)^n 
 \right] 
  =
   \sum_{ \rho \in \Pi [n\times r] } 
 \Bbb{E} \left[ 
 \int_{X^{|\rho |}} 
 \varepsilon^+_{\eufrak{z}_{|\rho|}} \prod_{k=1}^n
 u\big( z^\rho_{\pi_k} \big)
 \hat{\lambda}^{|\rho |} ( d \eufrak{z}_{|\rho|} ) \right] 
\end{equation} 
with 
$z^\rho_{\pi_k} := (z_{\zeta^\rho (k,1)}, \ldots , z_{\zeta^\rho (k,r)})$, 
$k=1,\ldots , n$. 
\end{prop} 
\begin{Proof}
  The main change in the proof argument of \cite{bogdan}
  is to rewrite the proof of Lemma~2.1 therein
  by applying \eqref{eq:hatc}
  recursively as in the proof of Theorem~3.1 of \cite{flint},
  while the main combinatorial argument remains identical.
  \end{Proof} 
  When $n=1$, Proposition~\ref{p2.2} yields a multivariate version of the
Georgii-Nguyen-Zessin identity
\eqref{eq:GNZ0}, i.e. 
$$ 
  \Bbb{E} \left[ 
  \int_{X^r} u(z_1,\ldots , z_r ; \omega ) \omega  (dz_1) \cdots
   \omega (dz_r) 
 \right] 
= 
 \sum_{\rho \in \Pi [1\times r] }
 \Bbb{E} \left[ 
 \int_{X^r} 
 \varepsilon^+_{\eufrak{z}_{|\rho|}} u(z_{\zeta^\rho (1,1)},\ldots ,z_{\zeta^\rho (1,r)} ; \omega )
 \hat{\lambda}^{|\rho |} ( d \eufrak{z}_{|\rho|} ) \right]. 
$$ 
\subsubsection*{Non-flat partitions} 
In the sequel we write $\nu \preceq \sigma$ when a partition 
$\nu \in \Pi [n \times r]$ is finer than another partition $\sigma \in \Pi [n \times r]$,
i.e. when every block of $\nu$ is contained in a block of $\sigma$.
We also let $\hat{0} : = \{ \{1,1\},\ldots , \{n,r\}\}$ denote the partition of
 $\Delta_{n\times r}$ made 
 of singletons, and we write $\rho \wedge \nu = \hat{0}$ when 
 $\mu = \hat{0}$ is the only partition
 $\mu \in \Pi [n \times r]$ such that
$\mu \preceq \nu$ and $\mu \preceq \rho$, i.e.
$|\nu_k \cap \rho_l| \leq 1$, $k=1,\ldots , n$,
$l=1, \ldots , |\rho|$.
 In this case the partition diagram of $\Gamma (\nu , \rho )$
 of $\nu$ and $\rho$ is said to be {\em non-flat},
 see Chapter~4 of \cite{peccatitaqqu}. 
 \\

 A partition $\rho$ is {\em non-flat} if the partition diagram of $\Gamma (\pi , \rho )$
 of $\rho$ and the partition  
 $\pi : = (\pi_1,\ldots , \pi_n) \in \Pi [n\times r]$
 with 
$\pi_k := \{ (k,1), \ldots , (k,r) \}$, $k=1,\ldots , n$, is {\em non-flat}.
 The following figure shows an example of a non-flat partition 

 
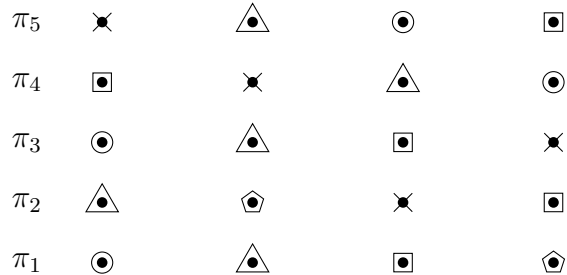
\begin{center}
   \begin{tikzpicture}[x=1cm,y=0.4cm]

\draw (-1,2) node {$\pi_1$};
\fill[black] (0,2) circle[radius=2pt];
\fill[black] (2,2) circle[radius=2pt];
\fill[black] (4,2) circle[radius=2pt];
\fill[black] (6,2) circle[radius=2pt];

\draw (0,2) circle (4pt);
\draw ([xshift=-3.5pt,yshift=-3.5pt]4,2) rectangle ++(7pt,7pt);
\node[draw,regular polygon, regular polygon sides=3,inner sep=2.5pt] at (2,2) {};
\node[draw,regular polygon, regular polygon sides=5,inner sep=2.5pt] at (6,2) {};

\draw (-1,4) node {$\pi_2$};
\fill[black] (0,4) circle[radius=2pt];
\fill[black] (2,4) circle[radius=2pt];
\fill[black] (4,4) circle[radius=2pt];
\fill[black] (6,4) circle[radius=2pt];

\draw ([xshift=-3.5pt,yshift=-3.5pt]6,4) rectangle ++(7pt,7pt);
\node[draw,regular polygon, regular polygon sides=3,inner sep=2.5pt] at (0,4) {};
\node[draw,regular polygon, regular polygon sides=5,inner sep=2.5pt] at (2,4) {};

\tikzset{cross/.style={cross out, draw, 
         minimum size=7pt, 
         inner sep=0pt, outer sep=0pt}}
\draw (4,4) node[cross,rotate=0] {};

\draw (-1,6) node {$\pi_3$};
\fill[black] (0,6) circle[radius=2pt];
\fill[black] (2,6) circle[radius=2pt];
\fill[black] (4,6) circle[radius=2pt];
\fill[black] (6,6) circle[radius=2pt];

\draw (0,6) circle (4pt);
\draw ([xshift=-3.5pt,yshift=-3.5pt]4,6) rectangle ++(7pt,7pt);
\node[draw,regular polygon, regular polygon sides=3,inner sep=2.5pt] at (2,6) {};
\draw (6,6) node[cross,rotate=0] {};

\draw (-1,8) node {$\pi_4$};
\fill[black] (0,8) circle[radius=2pt];
\fill[black] (2,8) circle[radius=2pt];
\fill[black] (4,8) circle[radius=2pt];
\fill[black] (6,8) circle[radius=2pt];

\draw (6,8) circle (4pt);
\draw ([xshift=-3.5pt,yshift=-3.5pt]0,8) rectangle ++(7pt,7pt);
\node[draw,regular polygon, regular polygon sides=3,inner sep=2.5pt] at (4,8) {};
\draw (2,8) node[cross,rotate=0] {};

\draw (-1,10) node {$\pi_5$};
\fill[black] (0,10) circle[radius=2pt];
\fill[black] (2,10) circle[radius=2pt];
\fill[black] (4,10) circle[radius=2pt];
\fill[black] (6,10) circle[radius=2pt];

\draw (4,10) circle (4pt);
\draw ([xshift=-3.5pt,yshift=-3.5pt]6,10) rectangle ++(7pt,7pt);
\node[draw,regular polygon, regular polygon sides=3,inner sep=2.5pt] at (2,10) {};
\draw (0,10) node[cross,rotate=0] {};

\end{tikzpicture}
 \end{center}
  
with $n=5$, $r=4$, and
 \begin{align*}
 & \triangle = \{ (1,2), (2,1) ,(2,2),(3,3),(4,2)\},
   \\
&   \ocircle = \{ (1,1), (3,1) ,(4,4),(5,3)\},
   \\
 & \Box = \{ (1,3), (2,4) ,(3,3),(4,1),(5,4)\},
   \\
 & \pentagon = \{ (1,4), (2,2)\},
   \\
   & \times = \{ (2,3), (3,4) ,(4,2),(5,1)\}
      \\
      & \pi_k = \{ (k,1), (k,2) ,(k,3),(k,4),(k,5)\},
      \quad
      k=1,2,3,4,5.
 \end{align*} 

\subsubsection*{Processes vanishing on diagonals}

 The next consequence of Proposition~\ref{p2.2} shows that 
 when $u(z_1,\ldots , z_r;\omega )$ vanishes on the diagonals in $X^r$,
 the moments of
 $\int_{X^r} u(z_1,\ldots , z_r;\omega ) \omega  (dz_1) \cdots
 \omega (dz_r)$
 reduce to sums over non-flat partition diagrams.  
 \begin{prop}
   \label{djks1} 
 Assume that $u(z_1,\ldots , z_r;\omega )=0$ whenever $z_i=z_j$,
 $1\leq i\not= j \leq r$, $\omega \in \Omega^X$. Then we have 
\begin{equation} 
\nonumber 
  \Bbb{E} \left[ 
   \left( \int_{X^r} u(z_1,\ldots , z_r;\omega ) \omega  (dz_1) \cdots
   \omega (dz_r) \right)^n 
 \right] 
= 
 \sum_{
   \substack{
     \rho \in \Pi [n\times r] 
     \\
     \rho \wedge \pi = \hat{0}
   }}
 \Bbb{E} \left[ 
 \int_{X^{|\rho |}} 
 \epsilon^+_{\eufrak{z}_{|\rho |}} \prod_{k=1}^n u \big( z^\rho_{\pi_k} \big) 
 \hat{\lambda}^{|\rho |} ( d \eufrak{z}_{|\rho |} ) \right]. 
\end{equation}
\end{prop} 
 When $n=1$, the first moment in Proposition~\ref{djks1}
 yields the Georgii-Nguyen-Zessin identity 
\begin{eqnarray} 
\nonumber
  \Bbb{E} \left[ 
   \int_{X^r} u(z_1,\ldots , z_r;\omega ) \omega  (dz_1) \cdots
   \omega (dz_r) 
 \right] 
& = & 
 \sum_{
   \substack{
     \rho \in \Pi [1\times r] 
     \\
     \rho \wedge \pi = \hat{0}
   }}
 \Bbb{E} \left[ 
 \int_{X^{|\rho |}} 
 \epsilon^+_{\eufrak{z}_{|\rho |}} u \big( z^\rho_{\pi_1} \big) 
 \hat{\lambda}^{|\rho |} ( d \eufrak{z}_{|\rho |} ) \right]
 \\
 \label{djksds} 
 & = & 
 \Bbb{E} \left[ 
 \int_{X^r} 
 \epsilon^+_{\eufrak{z}_r} 
 u(z_1,\ldots ,z_r ;\omega )
 \hat{\lambda}^r ( d \eufrak{z}_r ) \right]
,
\end{eqnarray} 
see Lemma~IV.1 in \cite{kartun-giles} and Lemma~2.1 in \cite{bogdan}
for different versions based on the Poisson point process.
In the case of second moments, we find 
$$ 
  \Bbb{E} \left[ 
   \left( \int_{X^r} u(z_1,\ldots , z_r;\omega ) \omega  (dz_1) \cdots
   \omega (dz_r) \right)^2 
 \right] 
 = 
 \sum_{
   \substack{
     \rho \in \Pi [2\times r] 
     \\
     \rho \wedge \pi = \hat{0}
   }}
 \Bbb{E} \left[ 
 \int_{X^{|\rho |}} 
 \epsilon^+_{\eufrak{z}_{|\rho |}} u \big( z^\rho_{\pi_1} \big) u \big( z^\rho_{\pi_2} \big) 
 \hat{\lambda}^{|\rho |} ( d \eufrak{z}_{|\rho |} ) \right] 
 ,
$$ 
 and since the non-flat partitions in $\Pi [ 2\times r]$
 are made of pairs and singletons, this identity
 can be rewritten as the following consequence of Proposition~\ref{djks1},
 in which for simplicity of notation we write $\pi_1 = \{1,\ldots ,r\}$ and
 $\pi_2 = \{r+1,\ldots ,2r\}$. 
\begin{corollary} 
\label{djklfs}
 Assume that $u(z_1,\ldots , z_r;\omega )=0$ whenever $z_i=z_j$,
 $1\leq i\not= j \leq r$, $\omega \in \Omega^X$. 
 Then the second moment of the integral of $k$-processes is given by 
\begin{align} 
\nonumber 
 & 
  \Bbb{E} \left[ 
   \left( \int_{X^r} u(z_1,\ldots , z_r;\omega ) \omega  (dz_1) \cdots
   \omega (dz_r) \right)^2 
 \right] 
\\
\nonumber
& = 
\sum_{
    A \subset \pi_1
}
\frac{1}{(r-|A|)!}
\sum_{
    \gamma : \pi_2 \to A \cup \{ r+|A|+1,\ldots , 2r\} 
}
\Bbb{E} \left[ 
 \int_{X^{2r-|A|}} \epsilon^+_{\eufrak{z}_{2r-|A|}} u \big( z_{\pi_1} \big) u \big( z_{\gamma (r+1)} , \ldots , z_{\gamma (2r)} \big) 
 \hat{\lambda}^{2r-|A|} ( d \eufrak{z}_{2r-|A|} ) \right]  
, 
\end{align} 
where the above sum is over all
bijections $\gamma : \pi_2 \to A \cup \{ r+|A|+1,\ldots , 2r\}$. 
\end{corollary}
\begin{Proof}
  We express the partitions $\rho \in \Pi [n\times r]$
  with non-flat diagrams $\Gamma ( \pi , \rho )$
  in Proposition~\ref{djklssa}
  as 
the collections of pairs and singletons
$$
\rho = \big\{ \{ i,\gamma (i) \}\}_{i\in A}
\cup 
\{\{i\}\}_{i\notin A}
\cup 
\{\{i\}\}_{i\notin \gamma (A)}
\big\},
$$ 
for all subsets $A\subset \pi_1 = \{1,\ldots , r \}$
and bijections $\gamma : \{1,\ldots , r\} \to A \cup \{ r+|A|+1,\ldots , 2r\}$.
\end{Proof}
 In the case of $2$-processes, Corollary~\ref{djklfs} shows that 
\begin{eqnarray} 
\nonumber 
\lefteqn{
  \Bbb{E} \left[ 
   \left( \int_{X^2} u(z_1 , z_2;\omega ) \omega  (dz_1) \omega (dz_2) \right)^2 
 \right] 
= 
  \sum_{
   \substack{
     \rho \in \Pi [n\times 2] 
     \\
     \rho \wedge \pi = \hat{0}
   }}
 \Bbb{E} \left[ 
 \int_{X^{|\rho |}} 
 \epsilon^+_{\eufrak{z}_{|\rho |}} \prod_{k=1}^n  u\big( z_{\zeta^\rho (k,1)},z_{\zeta^\rho (k,2)} \big)
 \hat{\lambda}^{|\rho |} ( d \eufrak{z}_{|\rho |} ) \right] 
}
\\
\nonumber 
& = &
\hskip-0.3cm 
 \sum_{
  \substack{
    A \subset \pi_1
    \\
    \gamma : \{ 3,4 \} \to A \cup \{ 3+|A|,\ldots , 4\}}
}
\hskip-0.3cm 
\frac{1}{(r-|A|)!}
\Bbb{E} \left[ 
 \int_{X^{4-|A|}} \epsilon^+_{\eufrak{z}_{4-|A|}} u ( z_1, z_2 ) u ( z_{\gamma (3)} , z_{\gamma (4)} ) 
 \hat{\lambda}^{4-|A|} ( d \eufrak{z}_{4-|A|} ) \right]  
\\
\nonumber
& = & 
 \Bbb{E} \left[ 
 \int_{X^4} 
 \epsilon^+_{\eufrak{z}_4} ( u(z_1,z_2 ) u(z_3,z_4 ) ) 
 \hat{\lambda}^4 ( d \eufrak{z}_4 ) \right]
\\
\nonumber
& & + 
 \Bbb{E} \left[ 
 \int_{X^3} 
 \epsilon^+_{\eufrak{z}_3} ( u(z_1,z_2 )
 u(z_1,z_3 ) ) 
 \hat{\lambda}^3 ( d \eufrak{z}_3 ) \right]
+ 
 \Bbb{E} \left[ 
 \int_{X^3} 
 \epsilon^+_{\eufrak{z}_3} ( u(z_2,z_1 )
 u(z_3,z_1 ) ) 
 \hat{\lambda}^3 ( d \eufrak{z}_3 ) \right]
\\
\nonumber
& &
+ 
 \Bbb{E} \left[ 
 \int_{X^3} 
 \epsilon^+_{\eufrak{z}_3} ( u(z_1,z_2 )
 u(z_2,z_3 ) ) 
 \hat{\lambda}^3 ( d \eufrak{z}_3 ) \right]
+ 
 \Bbb{E} \left[ 
 \int_{X^3} 
 \epsilon^+_{\eufrak{z}_3} ( u(z_2,z_1 )
 u(z_3,z_2 ) ) 
 \hat{\lambda}^3 ( d \eufrak{z}_3 ) \right]
 \\
\nonumber
& &
+ 
 \Bbb{E} \left[ 
 \int_{X^2} 
 \epsilon^+_{\eufrak{z}_2} ( u(z_1,z_2 )
 u(z_1,z_2 ) ) 
 \hat{\lambda}^2 ( d \eufrak{z}_2 ) \right]
 +
 \Bbb{E} \left[ 
 \int_{X^2} 
 \epsilon^+_{\eufrak{z}_2} ( u(z_1,z_2 )
 u(z_2,z_1 ) ) 
 \hat{\lambda}^2 ( d \eufrak{z}_2 ) \right]
. 
\end{eqnarray} 
 Similarly, in the case of $3$-processes we find 
\begin{align} 
\nonumber 
 & 
  \Bbb{E} \left[ 
    \left( \int_{X^3} u(z_1,z_1, z_3;\omega ) \omega  (dz_1) \omega (dz_2) \omega (dz_3)
    \right)^2 
 \right] 
\\
\nonumber
& = 
\sum_{
  \substack{
    A \subset \{1,2,3\}
    \\
    \gamma : \{4,5,6\} \to A \cup \{ 4+|A|,\ldots , 6\}}
}
\frac{1}{(3-|A|)!}
\Bbb{E} \left[ 
 \int_{X^5} \epsilon^+_{\epsilon^+_{\eufrak{z}_5} } u ( z_1,z_2,z_3 ) u \big( z_{\gamma (4)} , z_{\gamma (5)} ,z_{\gamma (6)} \big) 
 \hat{\lambda}^5 ( d \eufrak{z}_5  ) \right]  
\\
\nonumber
& = 
 \Bbb{E} \left[ 
 \int_{X^6} 
 \epsilon^+_{\eufrak{z}_6} u(z_1,z_2,z_3 ) u(z_4,z_5,z_6 )
 \hat{\lambda}^{6} ( d \eufrak{z}_6 ) \right]
 \\
\nonumber
& \quad + 
\frac{1}{2}
\sum_{ \gamma : \{4,5,6\} \to \{ 1,5 , 6\} }
\Bbb{E} \left[ 
 \int_{X^5} \epsilon^+_{\eufrak{z}_5} u ( z_1,z_2,z_3) u \big( z_{\gamma (4)} , z_{\gamma (5)} ,z_{\gamma (6)} \big) 
 \hat{\lambda}^5 ( d \eufrak{z}_5 ) \right]  
 \\
\nonumber
& \quad + 
\frac{1}{2}
\sum_{ \gamma : \{4,5,6\} \to \{ 2,5 , 6\} }
\Bbb{E} \left[ 
 \int_{X^5} \epsilon^+_{\eufrak{z}_5} u ( z_1,z_2,z_3) u \big( z_{\gamma (4)} , z_{\gamma (5)} ,z_{\gamma (6)} \big) 
 \hat{\lambda}^5 ( d \eufrak{z}_5 ) \right]  
 \\
\nonumber
& \quad + 
\frac{1}{2}
\sum_{ \gamma : \{4,5,6\} \to \{ 3,5 , 6\} }
\Bbb{E} \left[ 
 \int_{X^5} \epsilon^+_{\eufrak{z}_5} u ( z_1,z_2,z_3) u \big( z_{\gamma (4)} , z_{\gamma (5)} ,z_{\gamma (6)} \big) 
 \hat{\lambda}^5 ( d \eufrak{z}_5 ) \right]  
\\
\nonumber
& \quad + 
\sum_{
    \gamma : \{4,5,6\} \to \{1,2, 6\}
}
\Bbb{E} \left[ 
 \int_{X^4} \epsilon^+_{\eufrak{z}_4} u ( z_1,z_2,z_3 ) u \big( z_{\gamma (4)} , z_{\gamma (5)} ,z_{\gamma (6)} \big) 
 \hat{\lambda}^4 ( d \eufrak{z}_4 ) \right]  
\\
\nonumber
& \quad + 
\sum_{
    \gamma : \{4,5,6\} \to \{1,3, 6\}
}
\Bbb{E} \left[ 
 \int_{X^4} \epsilon^+_{\eufrak{z}_4} u ( z_1,z_2,z_3 ) u \big( z_{\gamma (4)} , z_{\gamma (5)} ,z_{\gamma (6)} \big) 
 \hat{\lambda}^4 ( d \eufrak{z}_4 ) \right]  
\\
\nonumber
& \quad + 
\sum_{
    \gamma : \{4,5,6\} \to \{2,3, 6\}
}
\Bbb{E} \left[ 
 \int_{X^4} \epsilon^+_{\eufrak{z}_4} u ( z_1,z_2,z_3 ) u \big( z_{\gamma (4)} , z_{\gamma (5)} ,z_{\gamma (6)} \big) 
 \hat{\lambda}^4 ( d \eufrak{z}_4 ) \right]  
\\
\nonumber
& \quad + 
\sum_{
 \gamma : \{4,5,6\} \to \{1,2,3 \} 
}
\Bbb{E} \left[ 
 \int_{X^3} \epsilon^+_{\eufrak{z}_3}z u ( z_1,z_2,z_3 ) u \big( z_{\gamma (4)} , z_{\gamma (5)} ,z_{\gamma (6)} \big) 
 \hat{\lambda}^3 ( d \eufrak{z}_3 ) \right]  
 . 
\end{align} 
\section{Random-connection model} 
\label{s4}
 Two point process vertices $x\not= y$ are connected
 in the random-connection graph with the probability  $H(x,y)$
 independently of $\omega (dx)$, 
 where $H:X\times X \longrightarrow [0,1]$. 
 In particular, the $1$-hop count $\mathds{1}_{\{ x \leftrightarrow y \} }$ is
 a Bernoulli random variable with parameter $H(x,y)$,
 and we have the relation  
$$ 
\Bbb{E} \left[ 
 \epsilon^+_{\eufrak{x}_r}
 \epsilon^+_{\eufrak{y}_r}
 \prod_{i=1}^r
 \prod_{j=1}^{i-1} 
 \mathds{1}_{\{ x_i \leftrightarrow y_j \} } ( \omega ) 
 \ \! \Big| \ \! \omega \right]
= 
\prod_{i=1}^r
\prod_{j=1}^{i-1} 
H( x_i , y_j ), 
$$
for distinct
$x_1,\ldots , x_r, y_1,\ldots , y_{r-1} \in X$,
where $x \leftrightarrow y$ means that $x\in X$ is connected to $y\in X$.
\\
 
Given $x,y \in X$, the number of $(r+1)$-hop sequences
$z_1,\ldots , z_r\in \omega$ of vertices
connecting $x$ to $y$ in the random graph is 
given by the multiparameter stochastic integral 
$$
 N^{x,y}_{r+1} = \int_{X^r} u(z_1,\ldots , z_r ;\omega ) \omega  (dz_1) \cdots \omega (dz_r)
$$
 of the $r$-process 
 \begin{equation}
   \label{fjkdsdf} 
u(z_1,\ldots , z_r;\omega )
: =
\mathds{1}_{\{
 z_i \not= z_j, \ 1\leq i<j \leq r \} }
\mathds{1}_{\{
  z_1,\ldots , z_r \in \omega \} }
\prod_{i=0}^r \mathds{1}_{\{ z_i \leftrightarrow z_{i+1} \} }
( \omega ),
\quad
z_1,\ldots , z_r \in X, 
\end{equation} 
which vanishes on the diagonals in $X^r$, with $z_0:=x$ and $z_{r+1} := y$.
In addition, for any distinct $z_1,\ldots , z_r \in X$ and
$u(z_1,\ldots , z_r;\omega )$ given by \eqref{fjkdsdf} we have 
\begin{equation}
   \label{djlsds} 
\Bbb{E} \big[ 
 \epsilon^+_{\eufrak{z}_r}
 u(z_1,\ldots , z_r;\omega ) 
 \mid \omega
 \big]
 = 
\Bbb{E} \left[ 
 \epsilon^+_{\eufrak{z}_r}
 \prod_{i=0}^r \mathds{1}_{\{ z_i \leftrightarrow z_{i+1} \} }
 ( \omega ) 
 \ \!
 \Big| \ \!
 \omega \right]
= 
\prod_{i=0}^r H( z_i , z_{i+1} )
, 
\end{equation} 
therefore the first order moment of the $(r+1)$-hop count
between $x\in X$ and $y\in X$ is given as  
\begin{equation}
  \label{jfdgfdg} 
 \Bbb{E} \left[ 
   \int_{X^r}
   u(z_1,\ldots , z_r; \omega ) 
    \omega  (dz_1) \cdots \omega (dz_r) 
 \right] 
 = 
 \Bbb{E} \left[
   \int_{X^r} 
 \prod_{i=0}^r H ( z_{1,i} , z_{1,i+1} ) 
 \hat{\lambda}^r ( d \eufrak{z}_r )
 \right],  
\end{equation} 
as a consequence of the Georgii-Nguyen-Zessin identity
 \eqref{djksds}. 
\\

 In the next proposition we compute the moments of all orders
of $r$-hop counts as sums over non-flat partition diagrams.
The role of the powers
$1/n^\rho_{l,i}$ in \eqref{dkasdad} 
is to ensure that all powers of $H(x,y)$ in \eqref{dkasdad} 
are equal to one, since all powers of 
$\mathds{1}_{\{ z \leftrightarrow z' \} }$ in \eqref{dhkjsd} below are equal
to $\mathds{1}_{\{ z \leftrightarrow z' \} }$. 
\begin{prop}
\label{djklssa} 
The moment of order $n$ of the $(r+1)$-hop count
between $x\in X$ and $y\in X$ is given by 
\begin{equation} 
\label{dkasdad} 
  \Bbb{E} \left[ 
    \left(
N^{x,y}_{r+1}  \right)^n 
 \right] 
 = 
 \sum_{
   \substack{
     \rho \in \Pi [n\times r] 
     \\
     \rho \wedge \pi = \hat{0}
   }}
 \Bbb{E} \left[
   \int_{X^{|\rho |}} 
 \prod_{l=1}^n
 \prod_{i=0}^r H^{1/n^\rho_{l,i}}( z_{\zeta^\rho (l,i)} , z_{\zeta^\rho (l,i+1)} ) 
 \hat{\lambda}^{|\rho |} ( d \eufrak{z}_{|\rho|} ) 
\right] ,
\end{equation} 
where $z_0=x$, $z_{r+1}=y$,
$\zeta^\rho (l,0)=0$, $\zeta^\rho (l,r+1)=r+1$, and
$$
 n^\rho_{l,i} := \# \big\{
 (p,j)\in \{1,\ldots ,n\}\times \{0,\ldots , r\} 
 \ : \
 \{ \zeta^\rho (l,i) , \zeta^\rho (l,i+1) \} 
 = \{ \zeta^\rho (p,j) , \zeta^\rho (p,j+1) \} 
 \big\}, 
 $$
 $1\leq l \leq n$,
 $0 \leq i \leq r$.
\end{prop}
\begin{Proof}
  Since $u(z_1,\ldots , z_r; \omega )$
  vanishes whenever $z_i=z_j$ for some $1\leq i < j \leq r$, 
  by Proposition~\ref{djks1} we have 
\begin{eqnarray} 
  \nonumber
  \lefteqn{
   \! \! \! \! \! \! \! \! \! \! \! \! \! \! \! \! \! \! \! \! \! \! \! \! \! \! 
  \Bbb{E} \left[ 
    \left(
    \int_{X^r}
    u(z_1,\ldots , z_r; \omega )
    \omega  (dz_1) \cdots \omega (dz_r) \right)^n 
 \right] 
  }
  \\
  \label{dhkjsd}
  & = &   
 \sum_{
   \substack{
     \rho \in \Pi [n\times r] 
     \\
     \rho \wedge \pi = \hat{0}
   }}
 \Bbb{E} \left[
   \int_{X^{|\rho |}} 
  \prod_{l=1}^n
 \prod_{i=0}^r \mathds{1}_{\{ z^\rho_{l,i} \leftrightarrow z^\rho_{l,i+1} \} }
 \hat{\lambda}^{|\rho |} ( d \eufrak{z}_{|\rho|} ) 
 \right]
\\
\nonumber
& = & 
 \sum_{
   \substack{
     \rho \in \Pi [n\times r] 
     \\
     \rho \wedge \pi = \hat{0}
   }}
 \Bbb{E} \left[
   \int_{X^{|\rho |}} 
  \prod_{l=1}^n
 \prod_{i=0}^r H^{1/n^\rho_{l,i}}( z^\rho_{l,i} , z^\rho_{l,i+1} ) 
 \hat{\lambda}^{|\rho |} ( d \eufrak{z}_{|\rho|} ) 
\right], 
\end{eqnarray} 
 where we applied \eqref{djlsds}.
\end{Proof} 
As in Corollary~\ref{djklfs} we have
the following consequence of Proposition~\ref{djklssa},
which is obtained by expressing the partitions $\rho \in \Pi [n\times r]$
with non-flat 
diagrams $\Gamma ( \pi , \sigma )$ as 
a collection of pairs and singletons.
\begin{corollary} 
\label{dvcsdf}
The second moment of the $(r+1)$-hop count
between $x\in X$ and $y\in X$ is given by 
\begin{eqnarray} 
\nonumber 
\lefteqn{
  \Bbb{E} \left[ 
    \left(
N^{x,y}_{r+1} \right)^2 
 \right] 
}
\\
\nonumber
& = &  
\hskip-1cm
\sum_{
  \substack{
    A \subset \pi_1
    \\
    \gamma : \{1,\ldots , r\} \to A \cup \{ r+|A|+1,\ldots , 2r\}}
}
 \hskip-1cm
 \frac{1}{(r-|A|)!}
 \Bbb{E} \left[
   \int_{X^{2r-|A|}}
 \prod_{i=0}^r H^{1/n^\gamma_{1,i}}( z_i , z_{i+1} ) 
 \prod_{j=0}^r H^{1/n^\gamma_{2,j}}( z_{\gamma (j)} , z_{\gamma (j+1) } ) 
 \hat{\lambda}^{2r-|A|} ( d \eufrak{z}_{2r-|A|} )
 \right], 
\end{eqnarray} 
where the above sum is over all bijections
$\gamma : \{1,\ldots , r\} \to A \cup \{ r+|A|+1,\ldots , 2r\}$ 
with $\gamma (0)=0$, $\gamma (r+1)=r+1$,
$z_0=x$, $z_{r+1}=y$,
and
 $$
 n^\gamma_{1,i} = \# \{
 j\in \{0,\ldots , r\} 
 \ : \
 \{ i , i+1 \} = \{ \gamma ( j) , \gamma (j+1) \}  
 \}, 
 $$
 $$
 n^\gamma_{2,j} = \# \{
 i\in \{0,\ldots , r\} 
 \ : \
 ( i , i+1 ) 
 = ( \gamma ( j) , \gamma (j+1) ) 
 \}, 
 $$
 $0 \leq i \leq r$.
\end{corollary}

\subsubsection*{Variance of $3$-hop counts}
When $n=2$ and $r=2$, 
Corollary~\ref{dvcsdf} allows us to
express the variance of the $3$-hop count 
between $x\in X$ and $y\in X$ as follows: 
\begin{align} 
\nonumber 
& 
 \Var \left[ 
        N^{x,y}_3 
 \right] 
\\
\nonumber
& =   
\hskip-0.8cm
\sum_{
  \substack{
 \emptyset \not= A \subset \{1,2\}
    \\
    \gamma : \{1,2\} \to A \cup \{ 3+|A| , 4\}}
}
 \hskip-0.8cm
 \frac{1}{(2-|A|)!}
 \Bbb{E} \left[
   \int_{X^{4-|A|}}
 \prod_{i=0}^2 H^{1/n^\gamma_{1,i}}( z_i , z_{i+1} ) 
 \prod_{j=0}^2 H^{1/n^\gamma_{2,j}}( z_{\gamma (j)} , z_{\gamma (j+1) } ) 
 \hat{\lambda}^{4-|A|} ( d \eufrak{z}_{4-|A|} ) 
 \right]
 \\
\nonumber
& =   
\sum_{
    \gamma : \{1,2\} \to \{ 1, 4\}
}
\Bbb{E} \left[
  \int_{X^3} 
 \prod_{i=0}^2 H^{1/n^\gamma_{1,i}}( z_i , z_{i+1} ) 
 \prod_{j=0}^2 H^{1/n^\gamma_{2,j}}( z_{\gamma (j)} , z_{\gamma (j+1) } ) 
 \hat{\lambda}^3 ( d z_1,dz_2,dz_4 ) 
 \right]
\\
\nonumber
&   
 + 
\sum_{
    \gamma : \{1,2\} \to \{ 2, 4\}
}
\Bbb{E} \left[
  \int_{X^3} 
 \prod_{i=0}^2 H^{1/n^\gamma_{1,i}}( z_i , z_{i+1} ) 
 \prod_{j=0}^2 H^{1/n^\gamma_{2,j}}( z_{\gamma (j)} , z_{\gamma (j+1) } ) 
 \hat{\lambda}^3 ( dz_1, dz_2, dz_4 ) 
 \right]
\\
\nonumber
&   
 + 
\sum_{
    \gamma : \{1,2\} \to \{ 1,2\}
}
\Bbb{E} \left[
  \int_{X^2} 
 \prod_{i=0}^2 H^{1/n^\gamma_{1,i}}( z_i , z_{i+1} ) 
 \prod_{j=0}^2 H^{1/n^\gamma_{2,j}}( z_{\gamma (j)} , z_{\gamma (j+1) } ) 
 \hat{\lambda}^2 ( dz_1,dz_2 ) 
 \right]
 .
\end{align} 

\subsubsection*{Variance of $4$-hop counts}
When $r=3$ and $n=2$, Corollary~\ref{dvcsdf} yields 
\begin{align*} 
 & 
\Var \left[ 
 N^{x,y}_4 
 \right] 
\\
\nonumber
& =   
\hskip-0.8cm
\sum_{
  \substack{
    \emptyset \not= A \subset \pi_1
    \\
    \gamma : \{1,\ldots , 3\} \to A \cup \{ 4+|A|,\ldots , 6\}}
}
 \hskip-0.8cm
 \frac{1}{(3-|A|)!}
 \Bbb{E} \left[
   \int_{X^{6-|A|}}
 \prod_{i=0}^3 H^{1/n^\gamma_{1,i}}( z_i , z_{i+1} ) 
 \prod_{j=0}^3 H^{1/n^\gamma_{2,j}}( z_{\gamma (j)} , z_{\gamma (j+1) } ) 
 \hat{\lambda}^{6-|A|} ( d \eufrak{z}_{6-|A|} )
 \right]
 \\
\nonumber
& =   
 \frac{1}{2} 
\sum_{
    \gamma : \{1,\ldots , 3\} \to \{ 1, 5, 6\}
}
\Bbb{E} \left[
  \int_{X^5} 
 \prod_{i=0}^3 H^{1/n^\gamma_{1,i}}( z_i , z_{i+1} ) 
 \prod_{j=0}^3 H^{1/n^\gamma_{2,j}}( z_{\gamma (j)} , z_{\gamma (j+1) } ) 
 \hat{\lambda}^5 ( dz_1,dz_2,dz_3,dz_5,dz_6 ) 
 \right]
\\
\nonumber
&  +   
 \frac{1}{2} 
\sum_{
    \gamma : \{1,\ldots , 3\} \to \{ 2, 5, 6\}
}
\Bbb{E} \left[
  \int_{X^5} 
 \prod_{i=0}^3 H^{1/n^\gamma_{1,i}}( z_i , z_{i+1} ) 
 \prod_{j=0}^3 H^{1/n^\gamma_{2,j}}( z_{\gamma (j)} , z_{\gamma (j+1) } ) 
 \hat{\lambda}^5 ( dz_1,dz_2,dz_3,dz_5,dz_6 ) 
 \right]
\\
\nonumber
&  +   
 \frac{1}{2} 
\sum_{
    \gamma : \{1,\ldots , 3\} \to \{ 3, 5, 6\}
}
\Bbb{E} \left[
  \int_{X^5} 
 \prod_{i=0}^3 H^{1/n^\gamma_{1,i}}( z_i , z_{i+1} ) 
 \prod_{j=0}^3 H^{1/n^\gamma_{2,j}}( z_{\gamma (j)} , z_{\gamma (j+1) } ) 
 \hat{\lambda}^5 ( dz_1,dz_2,dz_3,dz_5,dz_6 ) 
 \right]
\\
\nonumber
&  +   
\sum_{
    \gamma : \{1,\ldots , 3\} \to \{ 1,2, 6\} 
}
\Bbb{E} \left[
  \int_{X^4}
 \prod_{i=0}^3 H^{1/n^\gamma_{1,i}}( z_i , z_{i+1} ) 
 \prod_{j=0}^3 H^{1/n^\gamma_{2,j}}( z_{\gamma (j)} , z_{\gamma (j+1) } ) 
 \hat{\lambda}^4 ( dz_1,dz_2,dz_3,dz_6 ) 
 \right]
\\
 \nonumber
&  +   
\sum_{
    \gamma : \{1,\ldots , 3\} \to \{ 1,3, 6\} 
}
\Bbb{E} \left[
  \int_{X^4}
 \prod_{i=0}^3 H^{1/n^\gamma_{1,i}}( z_i , z_{i+1} ) 
 \prod_{j=0}^3 H^{1/n^\gamma_{2,j}}( z_{\gamma (j)} , z_{\gamma (j+1) } ) 
 \hat{\lambda}^4 ( dz_1,dz_2,dz_3,dz_6 ) 
 \right]
\\
\nonumber
&  +   
\sum_{
    \gamma : \{1,\ldots , 3\} \to \{ 2,3, 6\} 
}
\Bbb{E} \left[
  \int_{X^4}
 \prod_{i=0}^3 H^{1/n^\gamma_{1,i}}( z_i , z_{i+1} ) 
 \prod_{j=0}^3 H^{1/n^\gamma_{2,j}}( z_{\gamma (j)} , z_{\gamma (j+1) } ) 
 \hat{\lambda}^4 ( dz_1,dz_2,dz_3,dz_6 ) 
 \right]
\\
\nonumber
&  +   
\sum_{
    \gamma : \{1,\ldots , 3\} \to \{1,\ldots , 3\}
}
\Bbb{E} \left[
  \int_{X^3}
 \prod_{i=0}^3 H^{1/n^\gamma_{1,i}}( z_i , z_{i+1} ) 
 \prod_{j=0}^3 H^{1/n^\gamma_{2,j}}( z_{\gamma (j)} , z_{\gamma (j+1) } ) 
 \hat{\lambda}^3 ( dz_1,dz_2,dz_3 ) 
 \right]
. 
\end{align*} 
\section{Poisson case} 
\label{s4.1}
In this section and the next one, we 
work in the Poisson random-connection model,
using a Poisson point process on $X=\real^d$ with intensity
$\lambda (dx)$ on $\real^d$. 
We let 
\begin{equation}
  \label{djkldsf} 
 H^{(n)} (x_0,x_n ) : = \int_{\real^d} \cdots \int_{\real^d} 
 \prod_{i=0}^{n-1} H (x_i,x_{i+1})
 \lambda (dx_1) \cdots \lambda (dx_{n-1}), \quad
 x_0, x_n \in \real^d,
 \quad n \geq 1.  
\end{equation} 
 The $2$-hop count
 between $x\in \real^d$ and $y\in \real^d$
 is given by the first order stochastic integral 
$$ 
\int_{\real^d} u(z;\omega ) \omega  (dz) 
=
\int_{\real^d} \mathds{1}_{\{ x \leftrightarrow z_1 \} }
\mathds{1}_{\{ z_1 \leftrightarrow y \} }
( \omega ) \omega  (dz_1) 
=
\int_{\real^d} \mathds{1}_{\{ x \leftrightarrow z_1 \} }
\mathds{1}_{\{ z_1 \leftrightarrow y \} } \omega  (dz_1) 
,
$$ 
 and its moment of order $n$ is
\begin{eqnarray} 
\nonumber 
  \Bbb{E} \left[ 
   \left( \int_{\real^d} u(z_1;\omega ) \omega  (dz_1) \right)^n 
 \right] 
 & = &   
  \Bbb{E} \left[ 
    \left( \int_{\real^d} \mathds{1}_{\{ x \leftrightarrow z_1 \} }
    \mathds{1}_{\{ z_1 \leftrightarrow y \} }
    \omega  (dz_1) \right)^n 
 \right] 
\\
\nonumber
& = & 
   \sum_{
   \substack{
     \rho \in \Pi [n\times 1] 
   }}
 \int_{X^{|\rho |}} 
 \prod_{l=1}^{|\rho |}
 \big( H(x,z_l)H(z_l,y)
 \big) 
 \lambda^{|\rho |} ( dz_1,\ldots , dz_{|\rho|} ) 
\\
\nonumber
& = & 
\sum_{k=1}^n 
S(n,k)
 \left(
 \int_{\real^d} 
 H(x,z)H(z,y) 
 \lambda ( d z ) 
 \right)^k
\\
\nonumber
& = & 
\sum_{k=1}^n 
S(n,k)
\left(
H^{(2)} (x,y) 
 \right)^k
, 
\end{eqnarray} 
therefore, from \eqref{ws3},  
the $2$-hop count
between $x\in \real^d$ and $y\in \real^d$
is a Poisson random variable with mean
$$
H^{(2)}(x,y) = \int_{\real^d} H(x,z)H(z,y) \lambda ( d z ).
$$
 By \eqref{jfdgfdg}, the first order moment of the $r$-hop count is given by 
 $$
 H^{(r)}(x,y) = \int_{X^{r-1}} 
 \prod_{i=0}^{r-1} H ( z_i , z_{i+1} ) \lambda^{r-1} ( dz_1, \ldots dz_{r-1} ). 
$$ 
  
\begin{corollary} 
  \label{dvcsdf2}
  The variance of the $r$-hop count
  between $x\in \real^d$ and $y\in \real^d$ is given by 
\begin{align} 
\nonumber 
 & \Var \left[ 
N^{x,y}_r 
 \right] 
\\
\nonumber
& = 
\sum_{p=1}^{r-1}
\sum_{
  \substack{
    1\leq k_1  <\cdots < k_p < r 
  \\
  1 \leq  l_1  <\cdots < l_p < r
}
  }
\sum_{\sigma \in \Sigma [p]} 
 \int_{X^p} 
 \prod_{
      0 \leq i \leq p
 }
 H^{(k_{i+1}-k_i)}( z_i , z_{i+1} )
 \hskip-0.8cm
 \prod_{
   \substack{
     0 \leq j \leq p
   \\
   l_{\sigma (j+1)}-l_{\sigma (j)}
   + k_{j+1}-k_j > 2 
   \\
   \mbox{\scriptsize or } \{j,j+1\} \not= \{ \sigma (j), \sigma (j+1) \}}
 }
 \hskip-0.8cm
 H^{(l_{\sigma (j+1)}-l_{\sigma (j)})}( z_{\sigma (j)} , z_{\sigma ( j+1 ) } )
 \lambda^p ( d \eufrak{z}_p )
, 
\end{align} 
with $k_0=l_0=0$, $k_{p+1}=l_{p+1}=r$,
$\sigma (0)=0$, and $\sigma (r)=r$,
where the above sum if over all permutations
$\sigma \in \Sigma [p]$ of $\{1,\ldots , p\}$.
\end{corollary}
\begin{Proof}
  We rewrite the result of Corollary~\ref{dvcsdf} by
  denoting the set $A\subset \pi_1$ as
  $A = \{k_1, \ldots , k_p \}$,
  for $1\leq k_1  <\cdots < k_p \leq r-1$,
  and we identify $\gamma (A) \subset A \cup \{ r+|A|+1,\ldots , 2r\}$
  to $\{l_1, \ldots , l_p \}$,
  which requires a sum over the permutations of $\{1,\ldots , p\}$ since
  $1\leq l_1  <\cdots < l_p \leq r-1$,
  where $1 \leq p \leq r-1$.
  In addition, the multiple integrals over contiguous index sets
  in $A^c$ are evaluated using \eqref{djkldsf}. 
\end{Proof} 
\subsubsection*{Variance of $3$-hop counts}
When $n=2$ and $r=3$
Corollary~\ref{dvcsdf2} allows us to
compute the variance of the $3$-hop count
between $x\in \real^d$ and $y\in \real^d$, as follows: 
\begin{align} 
\label{djkdsjkfdsl} 
& 
 \Var \left[ 
        N^{x,y}_3 
 \right] 
\\
\nonumber
& =   
 2 \int_{\real^d} 
 H(x,z_1)H^{(2)}(z_1,y) H^{(2)}(z_1,y)
 \lambda ( d z_1 ) 
 + 
 2 \int_{\real^d}
 H(x,z_1)H^{(2)}(x,z_1)H^{(2)}(z_1,y)H(z_1,y)
 \lambda ( d z_1 ) 
\\
\nonumber
 & \quad   
 + \int_{X^2} H(x,z_1)H(z_1,z_2)H(z_2,y) H(x,z_2) H(z_1,y) \lambda^2 ( d z_1,dz_2 ) 
 + H^{(3)}(x,y)
 .
\end{align} 

\subsubsection*{Variance of $4$-hop counts}
 By Corollary~\ref{dvcsdf2} we have 
\begin{align} 
\label{fdjkldsfs} 
& 
\Var \left[ 
N^{x,y}_4 
 \right] 
\\
\nonumber
& =   
 \int_{\real^d} 
 H_\beta (x,z_1) H^{(3)}_\beta (z_1,y)
 H^{(3)}_\beta (z_1,y)
 \lambda ( d z_1 ) 
  +
 \int_{\real^d} 
 H_\beta (x,z_1)H^{(3)}_\beta (z_1,y)
 H^{(2)}_\beta (x,z_1)H^{(2)}_\beta (z_1,y)
 \lambda ( d z_1 ) 
 \\
\nonumber &  +
 \int_{\real^d}
 H_\beta (x,z_1)H^{(3)}_\beta (x,z_1)H^{(3)}_\beta (z_1,y)
 H_\beta (z_1,y)
 \lambda ( d z_1 ) 
 \\
\nonumber &  +
 \int_{\real^d} 
 H^{(2)}_\beta (x,z_2)H^{(2)}_\beta (z_2,y)
 H_\beta (x,z_2)H^{(3)}_\beta (z_2,y)
 \lambda ( d z_2 ) 
 \\
\nonumber &  +
 \int_{\real^d} 
 H^{(2)}_\beta (x,z_2)H^{(2)}_\beta (x,z_2)
 H^{(2)}_\beta (z_2,y)H^{(2)}_\beta (z_2,y)
 \lambda ( d z_2 ) 
 \\
\nonumber &  +
 \int_{\real^d} 
 H^{(2)}_\beta (x,z_2)H^{(2)}_\beta (z_2,y)H^{(3)}_\beta (x,z_2)H_\beta (z_2,y)
 \lambda ( d z_2 ) 
 \\
\nonumber &  +
 \int_{\real^d}
 H^{(3)}_\beta (x,z_3)H_\beta (z_3,y)H_\beta (x,z_3)H^{(3)}_\beta (z_3,y)
 \lambda ( d z_3 ) 
 \\
\nonumber &  +
 \int_{\real^d} 
 H^{(3)}_\beta (x,z_3)H_\beta (z_3,y)H^{(2)}_\beta (x,z_3)H^{(2)}_\beta (z_3,y)
 \lambda ( d z_3 ) 
 \\
\nonumber &  +
 \int_{\real^d} 
 H^{(3)}_\beta (x,z_3)H_\beta (z_3,y)H^{(3)}_\beta (x,z_3)
 \lambda ( d z_3 ) 
 \\
\nonumber  & + 
 \int_{X^2} 
 H_\beta (x,z_1)H_\beta (z_1,z_2)H^{(2)}_\beta (z_2,y)H^{(2)}_\beta (z_2,y)
 \lambda^2 ( dz_1,d z_2 ) 
 \\
\nonumber &  +
 \int_{X^2}
 H_\beta (x,z_1)H_\beta (z_1,z_2)H^{(2)}_\beta (z_2,y)H_\beta (x,z_2)H^{(2)}_\beta (z_1,y)
 \lambda^2 ( d z_1,dz_2 ) 
 \\
\nonumber &  +
 \int_{X^2} 
 H_\beta (x,z_1)H_\beta (z_1,z_2)H^{(2)}_\beta (z_2,y)H^{(2)}_\beta (z_1,z_2)H_\beta (z_2,y)
 \lambda^2 ( d z_1,dz_2 ) 
 \\
\nonumber &  +
 \int_{X^2} 
 H_\beta (x,z_1)H_\beta (z_1,z_2)H^{(2)}_\beta (z_2,y)H_\beta (x,z_2)H^{(2)}_\beta (z_2,z_1)H_\beta (z_1,y) 
 \lambda^2 ( d z_1,dz_2 ) 
 \\
\nonumber &  +
 \int_{X^2} 
 H_\beta (x,z_1)H_\beta (z_1,z_2)H^{(2)}_\beta (z_2,y)H^{(2)}_\beta (x,z_1)H_\beta (z_2,y)
 \lambda^2 ( d z_1,dz_2 ) 
 \\
\nonumber &  +
 \int_{X^2} 
 H_\beta (x,z_1)H_\beta (z_1,z_2)H_\beta (z_1,y)H^{(2)}_\beta (x,z_2)H^{(2)}_\beta (z_2,y)
 \lambda^2 ( d z_1,dz_2 ) 
 \\
\nonumber &  +
 \int_{X^2} 
 H_\beta (x,z_1)H^{(2)}_\beta (z_1,z_3)H_\beta (z_3,y)H_\beta (z_1,z_3)H^{(2)}_\beta (z_3,y)
 \lambda^2 ( d z_1,dz_3 ) 
 \\
\nonumber &  +
 \int_{X^2} 
 H_\beta (x,z_1)H^{(2)}_\beta (z_1,z_3)H^{(2)}_\beta (z_1,z_3)H_\beta (z_3,y)
 \lambda^2 ( d z_1,dz_3 ) 
 \\
\nonumber &  +
 \int_{X^2} 
 H_\beta (x,z_1)H^{(2)}_\beta (z_1,z_3)H_\beta (z_3,y)H_\beta (x,z_3)H_\beta (z_3,z_1)H^{(2)}_\beta (z_1,y)
 \lambda^2 ( d z_1,dz_3 ) 
 \\
\nonumber &  +
 \int_{X^2} 
 H_\beta (x,z_1)H^{(2)}_\beta (z_1,z_3)H_\beta (z_3,y)H_\beta (x,z_3)H^{(2)}_\beta (z_3,z_1)H_\beta (z_1,y)
 \lambda^2 ( d z_1,dz_3 ) 
 \\
\nonumber &  +
 \int_{X^2} 
 H_\beta (x,z_1)H^{(2)}_\beta (z_1,z_3)H^{(2)}_\beta (x,z_1)H_\beta (z_1,z_3)H_\beta (z_3,y)
 \lambda^2 ( d z_1,dz_3 ) 
 \\
\nonumber &  +
 \int_{X^2} 
 H_\beta (x,z_1)H^{(2)}_\beta (z_1,z_3)H_\beta (z_3,y)H^{(2)}_\beta (x,z_3)H_\beta (z_3,z_1)H_\beta (z_1,y)
 \lambda^2 ( d z_1,dz_3 ) 
 \\
\nonumber &  +
 \int_{X^2} 
 H^{(2)}_\beta (x,z_2)H_\beta (z_2,z_3)H_\beta (z_3,y)H_\beta (x,z_2)H^{(2)}_\beta (z_3,y)
 \lambda^2 ( d z_2,dz_3 ) 
 \\
\nonumber &  +
 \int_{X^2} 
 H^{(2)}_\beta (x,z_2)H_\beta (z_2,z_3)H_\beta (x,z_2)H^{(2)}_\beta (z_2,z_3)H_\beta (z_3,y)
 \lambda^2 ( d z_2,dz_3 ) 
 \\
\nonumber &  +
 \int_{X^2} 
 H^{(2)}_\beta (x,z_2)H_\beta (z_2,z_3)H_\beta (z_3,y)H_\beta (x,z_3)H^{(2)}_\beta (z_2,y)
 \lambda^2 ( d z_2,dz_3 ) 
 \\
\nonumber &  +
 \int_{X^2} 
 H^{(2)}_\beta (x,z_2)H_\beta (z_2,z_3)H_\beta (z_3,y)H_\beta (x,z_3)H^{(2)}_\beta (z_3,z_2)H_\beta (z_2,y)
 \lambda^2 ( d z_2,dz_3 ) 
 \\
\nonumber &  +
 \int_{X^2} 
 H^{(2)}_\beta (x,z_2)H_\beta (z_2,z_3)H_\beta (z_3,y)H^{(2)}_\beta (x,z_2)
 \lambda^2 ( d z_2,dz_3 ) 
 \\
\nonumber &  +
 \int_{X^2} 
 H^{(2)}_\beta (x,z_2)H_\beta (z_2,z_3)H_\beta (z_3,y)H^{(2)}_\beta (x,z_3)H_\beta (z_2,y)
 \lambda^2 ( d z_2,dz_3 ) 
 \\
\nonumber &  +
 H^{(4)}_\beta (x,y)
 + 
 \int_{X^3} 
 H_\beta (x,z_1)H_\beta (z_1,z_2)H_\beta (z_2,z_3)H_\beta(z_3,y)
 H_\beta (z_1,z_3)H_\beta(z_2,y)
 \lambda^3 ( d z_1,dz_2,dz_3 )
 \\
\nonumber &  +
 \int_{X^3} 
 H_\beta (x,z_1)H_\beta (z_1,z_2)H_\beta (z_2,z_3)H_\beta(z_3,y)
 H_\beta (x,z_2)H_\beta (z_1,z_3)
 \lambda^3 ( d z_1,dz_2,dz_3 )
 \\
\nonumber &  +
 \int_{X^3} 
 H_\beta (x,z_1)H_\beta (z_1,z_2)H_\beta (z_2,z_3)H_\beta(z_3,y)
 H_\beta (x,z_2)H_\beta (z_3,z_1)H_\beta(z_1,y)
 \lambda^3 ( d z_1,dz_2,dz_3 )
\\
\nonumber &  +
 \int_{X^3} 
 H_\beta (x,z_1)H_\beta (z_1,z_2)H_\beta (z_2,z_3)H_\beta(z_3,y)
 H_\beta (x,z_3)H_\beta (z_3,z_1)H_\beta(z_2,y)
 \lambda^3 ( d z_1,dz_2,dz_3 )
 \\
\nonumber &  +
 \int_{X^3} 
 H_\beta (x,z_1)H_\beta (z_1,z_2)H_\beta (z_2,z_3)H_\beta(z_3,y)
 H_\beta (x,z_3)H_\beta(z_1,y)
 \lambda^3 ( d z_1,dz_2,dz_3 )
. 
\end{align} 

\section{Rayleigh fading}
\label{s5}
In this section we consider a Poisson point process on $X=\real^d$ 
with flat intensity $\lambda (dx) = \lambda dx$ on $\real^d$,
$\lambda >0$, 
and a Rayleigh fading function of the form 
$$
H_\beta ( x , y ) := 
\re^{- \beta \Vert x - y \Vert^2 }, \qquad
x, y \in \real^d, 
\quad
\beta >0.
$$
Lemmas~\ref{djklsad}
and \ref{vcxc}
can be used to evaluate the
integrals appearing in Corollary~\ref{dvcsdf2} and in the variance
\eqref{fdjkldsfs} of $4$-hop counts. 
\begin{lemma} 
  \label{djklsad}
  For all $n\geq 1$, $y_1,\ldots, y_n \in \real^d$ and
  $\beta_1, \ldots , \beta_n >0$ we have 
  $$
  \int_{\real^d}
    \prod_{i=1}^n H_{\beta_i} (x,y_i)
    dx  
  = 
\left( \frac{\pi }{\beta_1+\cdots + \beta_n}\right)^{d/2}
\prod_{i=1}^{n-1} H_{\frac{\beta_{i+1} (\beta_1+\cdots + \beta_i)}{\beta_1+\cdots + \beta_{i+1}}} \left(
y_{i+1},
\frac{\beta_1 y_1+ \cdots + \beta_i y_i}{\beta_1+\cdots +\beta_i } \right)
.
$$ 
\end{lemma} 
\begin{Proof}
 We start by showing that for all $n\geq 1$ we have
 \begin{eqnarray} 
  \label{fdsfsfs}
 \lefteqn{
\prod_{i=1}^n H_{\beta_i} (x,y_i)
  }
\\
\nonumber
& = &
H_{\beta_1+\cdots + \beta_n} \left(x,\frac{\beta_1 y_1+ \cdots + \beta_n y_n}{\beta_1+\cdots +\beta_n }\right)
\prod_{i=1}^{n-1} H_{\frac{\beta_{i+1}(\beta_1+\cdots + \beta_i)}{\beta_1+\cdots + \beta_{i+1} }} \left(
y_{i+1},
\frac{\beta_1 y_1+ \cdots + \beta_i y_i}{\beta_1+\cdots +\beta_i } \right). 
\end{eqnarray} 
 Clearly, this relation holds for $n=1$.
 In addition, at the rank $n=2$ we have 
\begin{eqnarray*} 
  H_{\beta_1} (x,y_1)H_{\beta_2} (x,y_2) & = & 
   \re^{-\beta_1 \Vert y_1-x\Vert^2 } \re^{-\beta_2 \Vert x -y_2\Vert^2 }
\\
& = &
\re^{-\beta_1 \Vert y_1\Vert^2
  -\beta_2 \Vert y_2\Vert^2
  + 2 \langle \beta_1y_1 + \beta_2 y_2, x\rangle
  - ( \beta_1 + \beta_2 ) \Vert x\Vert^2 
}
\\
& = &
\re^{-\beta_1 \Vert y_1\Vert^2
  -\beta_2 \Vert y_2\Vert^2
  - ( \beta_1 + \beta_2 ) \Vert x - ( \beta_1y_1 + \beta_2 y_2)/(\beta_1+\beta_2 ) \Vert^2 
  + \Vert \beta_1y_1 + \beta_2 y_2 \Vert^2/(\beta_1+\beta_2 )  
}
\\
& = &
\re^{- ( \beta_1 + \beta_2 ) \Vert x - ( \beta_1y_1 + \beta_2 y_2)/(\beta_1+\beta_2 ) \Vert^2   - \beta_1 \beta_2 \Vert y_1-y_2\Vert^2 /(\beta_1+\beta_2 )  
}
\\
& = &
H_{\beta_1+\beta_2}
\left(x , \frac{\beta_1y_1 + \beta_2 y_2}{\beta_1+\beta_2} \right)
H_{\frac{ \beta_1 \beta_2 }{\beta_1+\beta_2 }}(y_1,y_2), 
\end{eqnarray*} 
 Next, assuming that \eqref{fdsfsfs} holds at the rank $n\geq 1$, we have 
    \begin{eqnarray*} 
  \lefteqn{
    \prod_{i=1}^{n+1} H_{\beta_i} (x,y_i)
    = H_{\beta_{n+1}} (x,y_{n+1})
H_{\beta_1+\cdots + \beta_n} \left(x,\frac{\beta_1 y_1+ \cdots + \beta_n y_n}{\beta_1+\cdots +\beta_n }\right)
}
\\
& &
\times \prod_{i=1}^{n-1} H_{\frac{\beta_{i+1}(\beta_1+\cdots + \beta_i)}{\beta_1+\cdots + \beta_{i+1}}} \left(
y_{i+1},
\frac{\beta_1 y_1+ \cdots + \beta_i y_i}{\beta_1+\cdots +\beta_i } \right)
\\
& = &
H_{\beta_1+\cdots + \beta_{n+1}} \left(x,\frac{\beta_1 y_1+ \cdots + \beta_{n+1} y_{n+1}}{\beta_1+\cdots +\beta_n }\right)
\prod_{i=1}^n H_{\frac{\beta_{i+1}(\beta_1+\cdots + \beta_i)}{\beta_1+\cdots + \beta_{i+1}}} \left(
y_{i+1},
\frac{\beta_1 y_1+ \cdots + \beta_i y_i}{\beta_1+\cdots +\beta_i } \right). 
\end{eqnarray*} 
 As a consequence, we find 
\begin{eqnarray*} 
  \lefteqn{
    \int_{\real^d}
    \prod_{i=1}^n H_{\beta_i} (x,y_i)
    dx
    =
  \prod_{i=1}^{n-1} H_{\frac{\beta_{i+1}(\beta_1+\cdots + \beta_i)}{\beta_1+\cdots + \beta_{i+1}}} \left(
y_{i+1},
\frac{\beta_1 y_1+ \cdots + \beta_i y_i}{\beta_1+\cdots +\beta_i } \right)
}
\\
& &
\times \int_{\real^d}
H_{\beta_1+\cdots + \beta_n} \left(x,\frac{\beta_1 y_1+ \cdots + \beta_n y_n}{\beta_1+\cdots +\beta_n }\right)
dx
\\
& = &
\left( \frac{\pi }{\beta_1+\cdots + \beta_n}\right)^{d/2}
\prod_{i=1}^{n-1} H_{\frac{\beta_{i+1}(\beta_1+\cdots + \beta_i)}{\beta_1+\cdots + \beta_{i+1}}} \left(
y_{i+1},
\frac{\beta_1 y_1+ \cdots + \beta_i y_i}{\beta_1+\cdots +\beta_i } \right)
.
\end{eqnarray*} 
\end{Proof}
  In particular, applying Lemma~\ref{djklsad} for $n=2$ yields 
\begin{eqnarray} 
  \label{djkljds}
  \int_{\real^d}
H_{\beta_1} ( y_1,x )
H_{\beta_2} ( x ,y_2 ) 
dx
& = &  
\left( \frac{\pi }{\beta_1 + \beta_2 }\right)^{d/2}
H_{\frac{\beta_1\beta_2}{\beta_1+\beta_2 }} ( y_1,y_2) 
\\
\nonumber
& = &  
\left( \frac{\pi }{\beta_1 + \beta_2 }\right)^{d/2}
e^{-  \beta_1\beta_2 \Vert y_1- y_2\Vert/ ( \beta_1+\beta_2 ) }, 
\quad
y_1,y_2 \in \real^d,
\end{eqnarray} 
and the $2$-hop count
between $x\in \real^d$ and $y\in \real^d$
is a Poisson random variable with mean
\begin{eqnarray} 
  \nonumber
  H^{(2)}_\beta (x,y) & = & \lambda \int_{\real^d}
H_\beta ( x,z )
H_\beta ( z ,y ) 
dz
\\
  \nonumber
& = &  
\lambda \left( \frac{\pi }{2 \beta }\right)^{d/2}
H_{\beta /2} ( x,y) 
\\
\nonumber
& = &  
\lambda \left( \frac{\pi }{2 \beta }\right)^{d/2}
e^{- \Vert x- y\Vert/ 2}. 
\end{eqnarray}
By an induction argument similar to that of
Lemma~\ref{djklsad}, we obtain the following lemma.
\begin{lemma} 
  \label{vcxc}
  For all $n \geq 1$, $x_0,\ldots, x_n \in \real^d$ and
  $\beta_1, \ldots , \beta_n >0$ we have 
\begin{eqnarray*} 
\lefteqn{
 \int_{\real^d} \cdots \int_{\real^d} 
 \prod_{i=1}^n H_{\beta_i} (x_{i-1},x_i)
 dx_1 \cdots dx_{n-1} 
}
\\
& = &
\left( \frac{\pi^{n-1} }{\sum_{i=1}^n \beta_1 \cdots \beta_{i-1}\beta_{i+1}\cdots \beta_n }\right)^{d/2}
H_{\frac{\beta_1 \cdots \beta_n}{\sum_{i=1}^n \beta_1 \cdots \beta_{i-1}\beta_{i+1}\cdots \beta_n }} (x_0,y_n). 
\end{eqnarray*} 
\end{lemma}
\begin{Proof}
  Clearly, the relation holds at the rank $n=1$.
  Assuming that it holds at the rank $n\geq 1$ and using \eqref{djkljds}, we have 
\begin{align*} 
 & 
 \int_{\real^d} \cdots \int_{\real^d} 
 \prod_{i=1}^{n+1} H_{\beta_i} (x_{i-1},x_i)
 dx_1 \cdots dx_n 
\\
& = 
\int_{\real^d}
H_{\beta_{n+1}} (x_n,x_{n+1})
\int_{\real^d} \cdots \int_{\real^d} 
 \prod_{i=1}^n H_{\beta_i} (x_{i-1},x_i)
 dx_1 \cdots dx_n 
 \\
 & = 
 \left( \frac{\pi^{n-1} }{\sum_{i=1}^n \beta_1 \cdots \beta_{i-1}\beta_{i+1}\cdots \beta_n }\right)^{d/2}
 \int_{\real^d}
 H_{\frac{\beta_1 \cdots \beta_n}{\sum_{i=1}^n \beta_1 \cdots \beta_{i-1}\beta_{i+1}\cdots \beta_n }} (x_0,x_n)
 H_{\beta_{n+1}} (x_n,x_{n+1})
 dx_n
 \\
 & = 
 \left( \frac{\pi^{n-1} }{\sum_{i=1}^n \beta_1 \cdots \beta_{i-1}\beta_{i+1}\cdots \beta_n }\right)^{d/2}
 \left( \frac{\pi }{\frac{\beta_1 \cdots \beta_n}{\sum_{i=1}^n \beta_1 \cdots \beta_{i-1}\beta_{i+1}\cdots \beta_n } + \beta_{n+1}}\right)^{d/2}
 H_{\frac{\beta_1 \cdots \beta_{n+1}}{\sum_{i=1}^{n+1} \beta_1 \cdots \beta_{i-1}\beta_{i+1}\cdots \beta_{n+1}}} ( x_0,x_{n+1}) 
.
\end{align*} 
\end{Proof} 
 In particular,
 the first order moment of the $r$-hop count between $x_0\in \real^d$ and $x_r\in \real^d$
 is given by 
\begin{eqnarray} 
\nonumber 
 H^{(r)}_\beta (x_0,x_r) & = & \int_{\real^d} \cdots \int_{\real^d} 
 \prod_{i=0}^{r-1} H_\beta (x_i,x_{i+1})
 \lambda ( dx_1 ) \cdots \lambda ( dx_{r-1} ) 
 \\
\nonumber 
 & = &  
\lambda^{r-1}
\left( \frac{\pi^{r-1} }{r \beta^{r-1} }\right)^{d/2}
H_{\beta / r } (x,y)
 \\
 \label{djks}
 & = &  
\lambda^{r-1}
\left( \frac{\pi^{r-1} }{r \beta^{r-1} }\right)^{d/2}
e^{-  \beta \Vert x- y\Vert/r},
\qquad x, y \in \real^d. 
\end{eqnarray} 
\subsubsection*{Variance of $3$-hop counts}
Corollary~\ref{dvcsdf2} and Lemma~\ref{vcxc} allow us to
recover Theorem~II.3 of \cite{kartun-giles}, 
for the variance of $3$-hop counts by a shorter argument, while extending it from
the plane $X = \real^2$ to $X = \real^d$.
\begin{corollary}
  \label{fdfdskl}
  The variance of the $3$-hop count
  between $x\in \real^d$ and $y\in \real^d$ is given by 
\begin{eqnarray*} 
    \Var \left[ 
    N^{x,y}_3 
 \right] 
&    = & 
    2 \lambda^3 \left( \frac{\pi^3}{8 \beta^3} \right)^{d/2}   
 e^{- \beta \Vert x- y\Vert/2}
 +
 \lambda^2 
\left( \frac{\pi^2 }{3 \beta^2 }\right)^{d/2}
 e^{- \beta \Vert x- y\Vert/3}
  \\
  & & +
 2 \lambda^3 \left( \frac{\pi^3}{12 \beta^3} \right)^{d/2}   
 e^{- 3 \beta \Vert x- y\Vert/4}
 + \lambda^2
 \left( \frac{\pi^2}{8 \beta^2 } \right)^{d/2} 
 e^{- \beta \Vert x- y\Vert }
. 
\end{eqnarray*} 
\end{corollary} 
\begin{Proof}
By \eqref{djks} and Lemma~\ref{vcxc} we have 
\begin{align*} 
& 
 \int_{\real^d} 
 H_\beta (x,z_1)H^{(2)}_\beta (z_1,y) H^{(2)}_\beta (z_1,y)
 \lambda ( d z_1 ) 
 = 
 \lambda^2
 \left( \frac{\pi^2}{4\beta^2} \right)^{d/2} 
 \int_{\real^d} 
 H_\beta (x,z_1)H^2_{\beta /2} (z_1,y) 
 \lambda ( d z_1 ) 
 \\
 & \quad
 = 
 \lambda^3
 \left( \frac{\pi^2}{4\beta^2} \right)^{d/2} 
 \int_{\real^d} 
 H_\beta (x,z_1)H_\beta (z_1,y) 
 \lambda ( d z_1 ) 
 = \lambda^3
 \left( \frac{\pi^3}{8\beta^3} \right)^{d/2}
 H_{\beta /2} (x,y) 
, 
\\
 & 
 \int_{\real^d} 
 H_\beta (x,z_1)H^{(2)}_\beta (x,z_1)H^{(2)}_\beta (z_1,y)H_\beta (z_1,y)
 \lambda ( d z_1 ) 
\\
\nonumber
& \quad = \lambda^2
 \left( \frac{\pi^2}{4\beta^2} \right)^{d/2} 
 \int_{\real^d} 
 H_{3\beta/2} (z_1,y) H_{3\beta /2} (x,z_1)
 \lambda ( d z_1 ) 
 = \lambda^3
\left( \frac{\pi^3}{12\beta^3} \right)^{d/2}
H_{3\beta /4} (x,y) 
, 
 \\
 & 
 \int_{X^2} H_\beta (x,z_1)H_\beta (z_1,z_2)H_\beta (z_2,y) H_\beta (x,z_2) H_\beta (z_1,y) \lambda^2 ( d z_1,dz_2 ) 
\\
 & \quad =   
\lambda \left( \frac{\pi}{3\beta} \right)^{d/2} 
H_\beta (x,y)
\int_{\real^d}
H_{2\beta/3 }(z_2,(x+y)/2)
H_{2\beta }(z_2,(x+y)/2)
\lambda ( d z_2 ) 
\\
 & \quad =   
\lambda^2
\left( \frac{\pi^2}{8 \beta^2 } \right)^{d/2}  
 H_\beta ( x , y ) 
,
\end{align*} 
and we conclude by \eqref{djkdsjkfdsl}.
\end{Proof} 

\footnotesize 

\def\cprime{$'$} \def\polhk#1{\setbox0=\hbox{#1}{\ooalign{\hidewidth
  \lower1.5ex\hbox{`}\hidewidth\crcr\unhbox0}}}
  \def\polhk#1{\setbox0=\hbox{#1}{\ooalign{\hidewidth
  \lower1.5ex\hbox{`}\hidewidth\crcr\unhbox0}}} \def\cprime{$'$}

\end{document}